\documentclass[12pt]{amsart}

\usepackage{amssymb, amscd, comment, txfonts}
\usepackage[all]{xy}
\usepackage{graphicx}

\numberwithin{equation}{section}

\newtheorem{thm}{Theorem}[section]

\newtheorem{lem}[thm]{Lemma}
\newtheorem{cor}[thm]{Corollary}

\theoremstyle{definition}
\newtheorem{defn}[thm]{Definition}

\theoremstyle{remark}
\newtheorem{rem}[thm]{Remark}
\newtheorem{exmp}[thm]{Example}

\renewcommand{\ker}{\operatorname{Ker}}
\renewcommand{\hom}{\operatorname{Hom}}

\newcommand{\N}{\mathbb{N}}
\newcommand{\Z}{\mathbb{Z}}
\newcommand{\Q}{\mathbb{Q}}
\newcommand{\R}{\mathbb{R}}
\newcommand{\C}{\mathbb{C}}
\newcommand{\F}{\mathbb{F}}

\DeclareMathOperator{\sgn}{sgn}

\DeclareMathOperator{\im}{Im}

\DeclareMathOperator{\tor}{Tor}

\DeclareMathOperator{\hdeg}{h-deg}
\DeclareMathOperator{\ldeg}{l-deg}

\begin{document}

%%%%%%% Title %%%%%%%%%%%%%%%%%%%%%%%%%%%%%%%%%%%%%%%%%%%%%%%%%%%%%%%
\title[Normalization of Twisted Alexander Invariants]
{Normalization of Twisted Alexander Invariants}
\author{Takahiro Kitayama}
\address{Department of Mathematics, Tokyo Institute of Technology,
2-12-1 Ookayama, Meguro-ku, Tokyo 152-8551, Japan}
\email{kitayama@math.titech.ac.jp}

\begin{abstract}
Twisted Alexander invariants of knots are well-defined up to multiplication of units.
We get rid of this multiplicative ambiguity via a combinatorial method and define normalized twisted Alexander invariants.
We then show that the invariants coincide with sign-determined Reidemeister torsion in a normalized setting, and refine the duality theorem.
We further obtain necessary conditions on the invariants for a knot to be fibered, and study behavior of the highest degrees of the invariants.
\end{abstract}

\maketitle

%%%%%%% Section 1 %%%%%%%%%%%%%%%%%%%%%%%%%%%%%%%%%%%%%%%%%%%%%%%%%%%
\section{Introduction}

Twisted Alexander invariants, which coincide with Reidemeister torsion (\cite{Ki}, \cite{KL}), were introduced for knots in the 3-sphere by Lin \cite{L} and generally for finitely presentable groups by Wada \cite{Wad}.
They were given a natural topological definition by using twisted homology groups in the notable work of Kirk and Livingston \cite{KL}.
Many properties of the classical Alexander polynomial $\Delta_K$ were subsequently extended to the twisted case and it was shown that the invariants have much information on the topological structure of a space.
For example, necessary conditions on twisted Alexander invariants for a knot to be fibered were given by Cha \cite{C}, Goda and Morifuji \cite{GM}, Goda, Kitano and Morifuji \cite{GKM}, and Friedl and Kim \cite{FK}.
Moreover, even sufficient conditions for a knot to be fibered were obtained by Friedl and Vidussi \cite{FV1, FV3}.

It is well known that $\Delta_K$ can be normalized, for instance, by considering the skein relation.
In this paper, we first obtain the corresponding result in twisted settings.
The twisted Alexander invariant $\Delta_{K, \rho}$ associated to a linear representation $\rho$ is well-defined up to multiplication of units in a Laurent polynomial ring.
We show that the ambiguity can be eliminated via a combinatorial method constructed by Wada and define the normalized twisted Alexander invariant $\widetilde{\Delta}_{K, \rho}$ (See Definition \ref{def_NTA} and Theorem \ref{thm_NTA}).

Turaev \cite{T2} defined sign-determined Reidemeister torsion by refining the sign ambiguity of Reidemeister torsion for an odd-dimensional manifold and showed that the other ambiguity depends on the choice of Euler structures.
We also normalize sign-determined Reidemeister torsion $T_{K, \rho}$ for a knot and define $\widetilde{T}_{K, \rho}(t)$.
Then we prove the equality
\[ \widetilde{\Delta}_{K, \rho}(t) = \widetilde{T}_{K, \rho}(t). \]
(See Theorem \ref{thm_NAR}.)
This shows that $\widetilde{\Delta}_{K, \rho}$ is a simple homotopy invariant and gives rise to a refined version of the duality theorem for twisted Alexander invariants.
(See Theorem \ref{thm_D}.)

As an application, we extend the above necessary conditions on $\widetilde{\Delta}_{K, \rho}$ for fibered knots.
We can define the highest degree and the coefficient of the highest degree term of $\widetilde{\Delta}_{K, \rho}$.
We show that these values are completely determined for fibered knots.
(See Theorem \ref{thm_F}.)
Finally, we obtain the following inequality which bounds the free genus $g_f(K)$ from below by the highest degree $\hdeg \widetilde{\Delta}_{K, \rho}$:
\begin{equation}
2 \hdeg \widetilde{\Delta}_{K, \rho} \leq n (2 g_f(K) - 1). \label{eq_F}
\end{equation}
(See Theorem \ref{thm_H}.)

This paper is organized as follows.
In the next section, we first review the definition of twisted Alexander invariants for knots.
We also describe how to compute them from a presentation of a knot group and the duality theorem for unitary representations.
In Section $3$, we review Turaev's sign-determined Reidemeister torsion and the relation with twisted Alexander invariants.
In Section $4$, we establish normalization of twisted Alexander invariants.
In Section $5$, we refine the correspondence with sign-determined Reidemeister torsion and the duality theorem for twisted Alexander invariants.
Section $6$ is devoted to applications.
Here we extend the result of Cha \cite{C}, Goda-Kitano-Morifuji \cite{GKM} and Friedl-Kim \cite{FK} for fibered knots, and study behavior of the highest degrees of the normalized invariants to obtain \eqref{eq_F}.

\subsection*{Note.}
This article appeared first in $2007$ on the arXiv, and has remained long to be unpublished.
Since then twisted Alexander invariants and Reidemeister torsion for knots and $3$-manifolds have been further intensively studied by many researchers.
We refer the reader to the survey papers~\cite{FV5, Mo} and the recent preprint~\cite{DFL} for details and references.
As this article has been already referred in the papers~\cite{DFJ, DFV, FKK, FV2, FV3, FV4, FV5, FV6, FV7, KM, SW} and frequently suggested to be published, we think that it might be worthwhile to have it  published. 

\subsection*{Acknowledgment.}
The author would like to express his gratitude to Toshitake Kohno for his encouragement and helpful suggestions.
The author would like to thank Hiroshi Goda, Teruaki Kitano, Takayuki Morifuji and Yoshikazu Yamaguchi for fruitful discussions, and Stefan Friedl for several stimulating comments which lead to some improvements of the argument in this revised version.
Finally the author also would like to thank the anonymous referee for helpful suggestions in revising the manuscript.
%The author was supported by JSPS KAKENHI (No.\ 26800032).
The author was supported by Research Fellowship of the Japan Society for the Promotion of Science for Young Scientists.

%%%%%%% Section 2 %%%%%%%%%%%%%%%%%%%%%%%%%%%%%%%%%%%%%%%%%%%%%%%%%%%
\section{Twisted Alexander invariants}

In this section, we review twisted Alexander invariants of an oriented knot, following \cite{C} and \cite{KL}.
For a given oriented knot $K$ in $S^3$, let $E_{K} := S^3 \setminus N(K)$, where $N(K)$ denotes an open tubular neighborhood of $K$, and let $G_K := \pi_1 E_K$.
We fix an element $\mu \in G_K$ represented by a meridian in $E_K$, and denote by $\alpha \colon G_K \to \langle t \rangle$ be the abelianization homomorphism which maps $\mu$ to the generator $t$.
Let $R$ be a Noetherian unique factorization domain and $Q(R)$ the quotient field of $R$.

We first define twisted homology groups and twisted cohomology groups.
Let $X$ be a connected CW-complex and $\widetilde{X}$ the universal cover of $X$.
The chain complex $C_*(\widetilde{X})$ is a left $\Z[\pi_1 X]$-module via the action of $\pi_1 X$ as deck transformations on $\widetilde{X}$.
We regard $C_*(\widetilde{X})$ also as a right $\Z[\pi_1 X]$-module by defining $\sigma \cdot \gamma := \gamma^{-1} \cdot \sigma$ for $\gamma \in \pi_1 X$ and $\sigma \in C_*(\widetilde{X})$.
For a linear representation $\rho \colon \pi_1 X \to GL_n(R)$, $R^{\oplus n}$ naturally has the structure of a left $\Z[\pi_1 X]$-module.
We define the \textit{twisted homology group} $H_i(X; R^{\oplus n}_{\rho})$ and the \textit{twisted cohomology group} $H^i(X; R^{\oplus n}_{\rho})$ associated to $\rho$ as follows:
\begin{align*}
H_i(X; R^{\oplus n}_{\rho}) &:= H_i(C_*(\widetilde{X}) \otimes_{\Z[\pi_1 X]} R^{\oplus n}), \\
H^i(X; R^{\oplus n}_{\rho}) &:= H^i(\hom_{\Z[\pi_1 X]}(C_*(\widetilde{X}), R^{\oplus n})).
\end{align*}

\begin{defn} \label{def_TA}
For a representation $\rho \colon G_K \to GL_n(R)$, we define $\Delta_{K, \rho}^{i}$ to be the order of the $i$-th twisted homology group $H_{i}(E_{K}; R[t,t^{-1}]_{\alpha \otimes \rho}^{\oplus n})$, where $R[t,t^{-1}]^{\oplus n} = R[t,t^{-1}] \otimes R^{\oplus n}$.
It is called the \textit{$i$-th twisted Alexander polynomial} associated to $\rho$, which is well-defined up to multiplication of units in $R[t,t^{-1}]$.
We furthermore define
\[ \Delta_{K, \rho} := \Delta_{K, \rho}^{1} / {\Delta_{K, \rho}^{0}} \in Q(R)(t), \]
which is called the \textit{twisted Alexander invariant} associated to $\rho$, and well-defined up to factors $\eta t^l$ for some $\eta \in R^{\times}$ and $l \in \Z$.
\end{defn}

\begin{rem}
Lin's twisted Alexander polynomial defined in \cite{L} coincides with $\Delta_{K, \rho}^{1}$.
\end{rem}

The homomorphisms $\alpha$ and $\alpha \otimes \rho$ naturally induce ring homomorphisms $\tilde{\alpha} \colon \Z [G_K] \to \Z [t,t^{-1}]$ and $\Phi \colon \Z [G_K] \to M_n(R[t,t^{-1}])$.
For a knot diagram of $K$, we choose and fix a Wirtinger presentation $G_K =  \langle x_1, \dots , x_m \mid r_1, \dots , r_{m-1} \rangle$.
Let us consider the $(m-1) \times m$ matrix $A_{\Phi}$ whose component is the $n \times n$ matrix $\Phi \left( \frac{\partial r_i}{\partial x_j} \right) \in M_n(R[t,t^{-1}])$, where $\frac{\partial}{\partial x_j}$ denotes Fox's free derivative with respect to $x_j$.
For $1 \leq k \leq m$, let us denote by $A_{\Phi, k}$ the $(m-1) \times (m-1)$ matrix obtained from $A_{\Phi}$ by removing the $k$-th column.
We naturally regard $A_{\Phi,k}$ as an $(m-1)n \times (m-1)n$ matrix with coefficients in $R[t,t^{-1}]$.

The twisted Alexander invariants can be computed from a Wirtinger presentation as follows.
The following is nothing but Wada's construction~\cite{Wad}.
\begin{thm}[\cite{HLN}, \cite{KL}] \label{thm_TA}
For a representation $\rho \colon G_K \to GL_n(R)$, a Wirtinger presentation $\langle x_1, \dots , x_m \mid r_1, \dots , r_{m-1} \rangle$ of $G_K$ and an index $k$,
\[ \Delta_{K,\rho} \equiv \frac{\det A_{\Phi,k}}{\det \Phi(x_k - 1)} \mod \langle \eta t^l \rangle_{\eta \in R^{\times}, l \in \Z}. \]
\end{thm}

\begin{rem}
Wada~\cite{Wad} showed that $\Delta_{K,\rho}$ is well-defined up to factors $\eta t^{ln}$.
He also showed that in the case where $\rho$ is a unimodular representation, $\Delta_{K,\rho}$ is well-defined up to factors $\pm t^{ln}$ if $n$ is odd and up to only $t^{ln}$ if $n$ is even. 
\end{rem}

It is also known that twisted Alexander invariants have the following duality.
We extend the complex conjugation to $\C(t)$ by taking $t \mapsto t^{-1}$.
\begin{thm}[\cite{Ki}, \cite{KL}] \label{thm_UD}
For a representation $\rho \colon G_K \to U(n)$ (resp. $O(n)$),
\[ \Delta_{K, \rho}(t) \equiv \overline{\Delta_{K, \rho}(t)} \mod \langle \eta t^l \rangle_{\eta \in R^{\times}, l \in \Z}. \]
\end{thm}

%%%%%%% Section 3 %%%%%%%%%%%%%%%%%%%%%%%%%%%%%%%%%%%%%%%%%%%%%%%%%%%
\section{Sign-determined Reidemeister torsion}

In this section, we review the definition of Turaev's sign-determined Reidemeister torsion.
See \cite{T1}, \cite{T2} for more details.
For two bases $u$ and $v$ of an $n$-dimensional vector space over a field $F$, $[u/v]$ denotes the determinant of the base change matrix from $v$ to $u$.

Let $C_* = (0 \to C_n \xrightarrow{\partial_n} C_{n-1} \to \cdots \xrightarrow{\partial_1} C_0 \to 0)$ be a chain complex of finite dimensional vector spaces over $F$.
For given bases $b_i$ of $\im \partial_{i+1}$ and $h_i$ of $H_i(C_{\ast})$, we can choose bases $b_i \cup \tilde{h}_i \cup \tilde{b}_{i-1}$ of $C_i$ as follows.
First, we choose a lift $\tilde{h}_i$ of $h_i$ in $\ker \partial_i$ and obtain a basis $b_i \cup \tilde{h}_i$ of $\ker \partial_i$, where we consider the exact sequence
\[ 0 \to \im \partial_{i+1} \to \ker \partial_i \to H_i(C_*) \to 0. \] 
Then we choose a lift $\tilde{b}_{i-1}$ of $b_{i-1}$ in $C_i$ and obtain a basis $(b_i \cup \tilde{h}_i) \cup \tilde{b}_{i-1}$ of $C_i$, where we consider the exact sequence
\[ 0 \to \ker \partial_i \to C_i \to \im \partial_i \to 0. \]

\begin{defn}
For given bases $\boldsymbol{c} = (c_i)$ of $C_*$ and $\boldsymbol{h} = (h_i)$ of $H_*(C_*)$, we choose a basis $\boldsymbol{b} = (b_i)$ of $\im \partial_*$ and define
\[ \tor (C_*,\boldsymbol{c},\boldsymbol{h}) := (-1)^{|C_*|} \prod_{i=0}^n [b_i \cup \tilde{h}_i \cup \tilde{b}_{i-1} / c_i]^{(-1)^{i+1}} ~ \in F^{\times}, \]
where
\[ |C_*| := \sum_{j=0}^n (\sum_{i=0}^j \dim C_i)(\sum_{i=0}^j \dim H_i(C_*)). \]
\end{defn}

\begin{rem}
It can be easily checked that $\tor (C_*,\boldsymbol{c},\boldsymbol{h})$ does not depend on the choices of $\boldsymbol{b}$, $\tilde{b}_i$ and $\tilde{h}_i$.
\end{rem}

Now let us apply the above algebraic torsion to geometric situations.
Let $X$ be a connected finite CW-complex.
By a \textit{homology orientation} of $X$ we mean an orientation of the homology group $H_*(X; \R) = \bigoplus_i H_i(X; \R)$ as a real vector space.

\begin{defn}
For a representation $\rho \colon \pi_1 X \to GL_n(F)$ such that $H_*(X; F_{\rho}^{\oplus n})$ vanishes and a homology orientation $\mathfrak{o}$, we define the \textit{sign-determined Reidemeister torsion} $T_{\rho}(X, \mathfrak{o})$ associated to $\rho$ and $\mathfrak{o}$ as follows.
We choose a lift $\tilde{e}_i$ of each cell $e_i$ in $\widetilde{X}$ and bases $\boldsymbol{h}$ of $H_*(X; \R)$ which is positively oriented with respect to $\mathfrak{o}$ and $\langle f_1, \dots, f_n \rangle$ of $F^{\oplus n}$.
Then we define
\[ T_{\rho}(X, \mathfrak{o}) := \tau_0^n \tor (C_*(\widetilde{X}) \otimes_{\rho} F^{\oplus n}, \tilde{\boldsymbol{c}}) ~ \in F^{\times}, \]
where
\begin{align*}
\tau_0 &:= \sgn \tor (C_*(X;\R), \boldsymbol{c}, \boldsymbol{h}), \\
\boldsymbol{c} &:= \langle e_1, \dots , e_{dim C_{\ast}} \rangle, \\
\tilde{\boldsymbol{c}} &:= \langle \tilde{e}_1 \otimes f_1, \dots , \tilde{e}_1 \otimes f_n, \dots , \tilde{e}_{dim C_{\ast}} \otimes f_1, \dots , \tilde{e}_{dim C_{\ast}} \otimes f_n \rangle.
\end{align*}
\end{defn}

\begin{rem}
It is known that $T_{\rho}(X, \mathfrak{o})$ does not depend on the choices of $\tilde{e}_i$, $\boldsymbol{h}$ and $\langle f_1, \dots, f_n \rangle$ and is well-defined as a simple homotopy invariant up to multiplication of elements in $\im(\det \circ \rho)$.
\end{rem}

Here let us consider the knot exterior $E_K$.
In this case, we can equip $E_K$ with its \textit{canonical homology orientation} $\omega_K$ as follows.
We have $H_*(E_K;\R) = H_0(E_K;\R) \oplus \langle t \rangle$, and define $\omega_K := [ \langle [pt], t \rangle ]$, where $[pt]$ is the homology class of a point.

\begin{defn}
For a representation $\rho \colon G_K \to GL_n(F)$ such that $H_*(X; F(t)_{\alpha \otimes \rho}^{\oplus n})$ vanishes, the \textit{sign-determined Reidemeister torsion} $T_{K,\rho}(t)$ associated to $\rho$ is defined by $T_{\alpha \otimes \rho}(E_K,\omega_K)$.
Here we consider $\alpha \otimes \rho$ as a representation $G_K \to GL_n(F[t,t^{-1}]) \hookrightarrow GL_n(F(t))$.
\end{defn}

In Section $5$, we generalize the following theorem.
\begin{thm}[\cite{Ki}, \cite{KL}] \label{thm_AR}
For a representation $\rho \colon G_K \to GL_n(F)$ such that $H_*(X; F(t)_{\alpha \otimes \rho}^{\oplus n})$ vanishes,
\[ \Delta_{K, \rho}(t) \equiv T_{K, \rho}(t) \mod \langle \eta t^l \rangle_{\eta \in F^{\times}, l \in \Z}. \]
\end{thm}

%%%%%%% Section 4 %%%%%%%%%%%%%%%%%%%%%%%%%%%%%%%%%%%%%%%%%%%%%%%%%%%
\section{Construction}

Now we establish one of our main results.
We get rid of the multiplicative ambiguity of twisted Alexander invariants via a combinatorial method.
For $f(t) = p(t) / q(t) \in Q(R)(t)$ $(p,q \in R[t,t^{-1}])$, we define
\begin{align*}
\deg f &:= \deg p - \deg q, \\
\hdeg f &:= (\text{the highest degree of } p) - (\text{the highest degree of } q), \\
\ldeg f &:= (\text{the lowest degree of } p) - (\text{the lowest degree of } q), \\
\mathrm{c}(f) &:= \frac{(\text{the coefficient of the highest degree term of } p)}{(\text{the coefficient of the highest degree term of } q)}.
\end{align*}

We make use of a combinatorial group theoretical approach constructed by Wada \cite{Wad}.
\begin{defn} \label{def_STT}
For a finite presentable group $G = \langle x_1, \dots , x_m~|~r_1, \dots , r_n \rangle $ and any word $w$ in $x_1, \dots , x_m$, the operations of the following types are called the \textit{strong Tietze transformations}:
\begin{itemize}
\item[Ia.] To replace one of the relators $r_i$ by its inverse $r_i^{-1}$.
\item[Ib.] To replace one of the relators $r_i$ by its conjugate $wr_iw^{-1}$.
\item[Ic.] To replace one of the relators $r_i$ by $r_ir_j$ for any $j \neq i$.
\item[II.] To add a new generator $y$ and a new relator $yw^{-1}$. (Namely, the resulting presentation is $\langle x_1, \dots , x_m, y~|~r_1, \dots ,r_n , yw^{-1} \rangle$.)
\end{itemize}
If one presentation is transformable to another by a finite sequence of operations of above types and their inverse operations, then such two presentations are said to be \textit{strongly Tietze equivalent}.
\end{defn}

\begin{rem}
The deficiency of a presentation does not change via the strong Tietze transformations.
\end{rem}

Wada showed the following lemma.
\begin{lem}[\cite{Wad}] \label{lem_W}
All the Wirtinger presentations for a given link in $S^3$ are strongly Tietze equivalent to each other.
\end{lem}

Let $\varphi \colon \Z [G_K] \to \Z$ be the augmentation homomorphism, namely, $\varphi (\gamma) = 1$ for any element $\gamma$ of $G_K$.
For a given presentation $\langle x_1, \dots , x_m~ |$ $~r_1, \dots , r_{m-1} \rangle$ of $G_K$, we denote $A_{\varphi,k}$ and $A_{\tilde{\alpha},k}$ by $\left( \varphi \left( \frac{\partial r_i}{\partial x_j} \right) \right)_{j \neq k}$ and $\left( \tilde{\alpha} \left( \frac{\partial r_i}{\partial x_j} \right) \right)_{j \neq k}$ as in Section $2$.

We eliminate the ambiguity of $\eta t^l$ in Definition \ref{def_TA} as follows.
\begin{defn} \label{def_NTA}
Given a representation $\rho \colon G_K \to GL_n(R)$, we choose a presentation $\langle x_1, \dots , x_m~|~r_1, \dots , r_{m-1} \rangle$ of $G_K$ which is strongly Tietze equivalent to a Wirtinger presentation and an index $k$ such that $\hdeg \alpha (x_k) \neq 0$.
Then we define the \textit{normalized twisted Alexander invariant} associated to $\rho$ as:
\[ \widetilde{\Delta}_{K,\rho} := \frac{\delta^n}{(\epsilon t^n)^d} \frac{\det A_{\Phi,k}}{\det \Phi(x_k - 1)} ~ \in Q(R)(\epsilon^{\frac{1}{2}})(t^{\frac{1}{2}}), \]
where
\begin{align*}
\epsilon &:= \det \rho (\mu), \\
\delta &:= \sgn (\hdeg \alpha (x_k) \det A_{\varphi,k}), \\
d &:= \frac{1}{2}(\hdeg \det A_{\tilde{\alpha},k} + \ldeg \det A_{\tilde{\alpha},k} - \hdeg \alpha (x_k)).
\end{align*}
\end{defn}

\begin{thm} \label{thm_NTA}
$\widetilde{\Delta}_{K,\rho}$ is an invariant of a linear representation $\rho$. 
\end{thm}

\begin{proof}
From Lemma \ref{lem_W}, it suffices to check $(i)$ the independence of the choice of $k$ and $(ii)$ the invariance for each operation in Definition \ref{def_STT}.

We assume that there is another index, say $k'$, also satisfying the condition $\hdeg \alpha (x_{k'}) \neq 0$.
We set
\begin{align*}
\delta' &:= \sgn (\hdeg \alpha (x_{k'}) \det A_{\varphi, {k'}}), \\
d' &:= \frac{1}{2}(\hdeg \det A_{\tilde{\alpha}, k'} + \ldeg \det A_{\tilde{\alpha}, k'} - \hdeg \alpha (x_{k'})).
\end{align*}
Since
\[ \sum\limits_{j=1}^m \frac{\partial r_i}{\partial x_j} (x_j - 1) = r_i - 1, \]
we have
\begin{align*}
\det A_{\Phi, k'} \det \Phi(x_k - 1) &= \det \left( \dots , \Phi \left( \frac{\partial r_i}{\partial x_k} \right) \Phi(x_k - 1) , \dots \right), \\
&= \det \left( \dots , - \sum_{j \neq k} \Phi \left( \frac{\partial r_i}{\partial x_j} \right) \Phi(x_j - 1) , \dots \right), \\
&= \det \left( \dots , -\Phi \left( \frac{\partial r_i}{\partial x_{k'}} \right) \Phi(x_{k'} - 1) , \dots \right), \\
&= (-1)^{n(k-k')} \det A_{\Phi,k} \det \Phi(x_{k'} - 1).
\end{align*}
Similarly, we have
\[ \det A_{\tilde{\alpha}, k'} \det \tilde{\alpha}(x_k - 1) = (-1)^{k - k'} \det A_{\tilde{\alpha}, k} \det \tilde{\alpha}(x_{k'} - 1). \]
Hence $d' = d$.
Moreover, by dividing this equality by $(t-1)$ and taking $t \to 1$, we can see that
\[ \hdeg \alpha(x_k) \det A_{\varphi, k'} = (-1)^{k-k'} \hdeg \alpha(x_{k'}) \det A_{\varphi, k}. \]
Hence $\delta' = (-1)^{k-k'} \delta$.
The above equalities prove $(i)$.

Next, we consider the strong Tietze transformations.
Since
\begin{align*}
\frac{\partial (r_i^{-1})}{\partial x_j} &= -r_i\frac{\partial r_i}{\partial x_j}, \\
\frac{\partial (wr_iw^{-1})}{\partial x_j} &= w\frac{\partial r_i}{\partial x_j}, \\
\frac{\partial (r_ir_l)}{\partial x_j} &= \frac{\partial r_i}{\partial x_j} + r_i\frac{\partial r_l}{\partial x_j},
\end{align*}
the changes of the values $\det A_{\Phi, k}$, $\delta$, $d$ by the transformation Ia, Ib and Ic are as follows.
By the transformation Ia, $\det A_{\Phi, k} \mapsto (-1)^n \det A_{\Phi, k}$, $\delta \mapsto -\delta$ and $d$ does not change.
By the transformation Ib, $\det A_{\Phi, k} \mapsto (\epsilon t^n)^{\deg \alpha(w)} \det A_{\Phi, k}$, $\delta$ does not change and $d \mapsto d + \deg \alpha(w)$.
By the transformation Ic and II, it is easy to see that all the values do not change.
These observations proves $(ii)$.
\end{proof}

From the construction, the following lemma holds.
\begin{lem} \label{lem_NA}
(i)For a representation $\rho \colon G_K \to GL_n(R)$,
\[ \Delta_{K, \rho}(t) \equiv \widetilde{\Delta}_{K, \rho}(t) \bmod \langle \epsilon^{\frac{1}{2}}, \eta t^{\frac{l}{2}} \rangle_{\eta \in R^{\times}, l \in \Z }. \]
(ii)If $\rho$ is trivial (i.e., $\Phi = \tilde{\alpha}$), then
\[ \nabla_K(t^{\frac{1}{2}}-t^{-\frac{1}{2}}) = (t^{\frac{1}{2}} - t^{-\frac{1}{2}}) \widetilde{\Delta}_{K, \rho}(t) ,\]
where $\nabla_K(z)$ is the Conway polynomial of $K$.
\end{lem}

\begin{proof}
Since (i) is clear from Theorem \ref{thm_TA} and Definition \ref{def_NTA}, we prove (ii).
For the trivial representation $\rho$, we set
\[ f(t) = (t^{\frac{1}{2}} - t^{-\frac{1}{2}}) \widetilde{\Delta}_{K, \rho}(t). \]
Then it is easy to see that
\[ f(t) \equiv \Delta_K(t) \mod \langle \pm t \rangle. \]
Moreover, we can check the following:
\begin{align*}
f(1) &= 1, \\
\hdeg f + \ldeg f &= 0,
\end{align*}
which establishes the desired formula.
\end{proof}

%%%%%%% Section 5 %%%%%%%%%%%%%%%%%%%%%%%%%%%%%%%%%%%%%%%%%%%%%%%%%%%
\section{Relation to sign-determined Reidemeister torsion}

In this section, we generalize Theorem \ref{thm_UD} and Theorem \ref{thm_AR}.
Here we only consider the case where $R$ is a field $F$.

First, we also normalize sign-determined Reidemeister torsion as twisted Alexander invariants.
\begin{defn} \label{def_NR}
For a representation $\rho \colon G_K \to GL_n(F)$ such that $H_*(E_K; F(t)_{\alpha \otimes \rho}^{\oplus n})$ vanishes, we define $\widetilde{T}_{K,\rho}(t)$ as follows.
We choose a lift $\tilde{e}_i$ in $\widetilde{E}_K$ of each cell $e_i$, bases $\boldsymbol{h}$ of $H_*(E_K;\R )$ which is positively oriented with respect to $\omega_K$ and $\langle f_1, \dots, f_n \rangle$ of $F(t)^{\oplus n}$.
Then we define
\[ \widetilde{T}_{K,\rho}(t) := \frac{\tau_0^n}{(\epsilon t^n)^{d'}} \tor (C_*(\widetilde{E}_K) \otimes_{\alpha \otimes \rho} F(t)^{\oplus n}, \tilde{\boldsymbol{c}} ) ~ \in F(t)^{\times}, \]
where
\begin{align*}
\epsilon &:= \det \rho (\mu), \\
\tau_0 &:= \sgn \tor (C_*(E_K;\R),\boldsymbol{c}, \boldsymbol{h}), \\
d' &:= \frac{1}{2} (\hdeg \tor (C_*(\widetilde{E_K}) \otimes_{\alpha} \Q (t), \tilde{\boldsymbol{c}}_0) + \ldeg \tor (C_*(\widetilde{E_K}) \otimes_{\alpha} \Q (t), \tilde{\boldsymbol{c}}_0)), \\
\boldsymbol{c} &:= \langle e_1, \dots , e_{dim C_*} \rangle, \\
\tilde{\boldsymbol{c}}_0 &:= \langle \tilde{e}_1 \otimes 1, \dots , \tilde{e}_{dim C_*} \otimes 1 \rangle, \\
\tilde{\boldsymbol{c}} &:= \langle \tilde{e}_1 \otimes f_1, \dots , \tilde{e}_1 \otimes f_n, \dots , \tilde{e}_{dim C_*} \otimes f_1, \dots , \tilde{e}_{dim C_*} \otimes f_n \rangle.
\end{align*}
\end{defn}

\begin{rem}
We can also define normalized Reidemeister torsion for an oriented link whose Alexander polynomial does not vanish by a similar method as follows:
When $K$ is an oriented link with ordered components $K_1, \dots, K_m$, we think $\alpha \colon G_K \to \langle t_1, \dots, t_m \rangle$ as the homomorphism which maps the meridional element $\mu_i$ of $K_i$ to the generator $t_i$ for each $i$, and define the canonical homology orientation as $\omega_K := [\langle [pt], [\mu_1], \dots, [\mu_m] \rangle]$. 
In the notation in Definition \ref{def_NR} we replace the field $F(t)$ by $F(t_1, \dots, t_m)$, and instead of $\epsilon$ and $d$ we set
\begin{align*}
\epsilon_i &:= \det \rho (\mu_i), \\
d_i' &:= \frac{1}{2} (\hdeg_i \tor (C_*(\widetilde{E_K}) \otimes_{\alpha} \Q (t_1, \dots, t_m), \tilde{\boldsymbol{c}}_0) + \ldeg_i \tor (C_*(\widetilde{E_K}) \otimes_{\alpha} \Q (t_1, \dots, t_m), \tilde{\boldsymbol{c}}_0)),
\end{align*}
where $\hdeg_i$ and $\ldeg_i$ are defined as $\hdeg$ and $\ldeg$ for polynomials on $t_i$.
Then we define
\[ \widetilde{T}_{K,\rho}(t_1, \dots, t_m) := \frac{\tau_0^n}{(\epsilon_1 t_1^n)^{d_1'} \cdots (\epsilon_m t_m^n)^{d_m'}} \tor (C_*(\widetilde{E}_K) \otimes_{\alpha \otimes \rho} F(t_1, \dots, t_m)^{\oplus n}, \tilde{\boldsymbol{c}} ) ~ \in F(t_1, \dots, t_m)^{\times}. \]
Note that if we permute two of the indices of components $K_1, \dots, K_m$, then the normalized invariant is multiplied with $(-1)^n$.
\end{rem}

One can prove the following lemma by a similar way as in the non-normalized case.
As a reference, see \cite{T1}.
\begin{lem} \label{lem_NR}
$\widetilde{T}_{K,\rho}$ is invariant under homology orientation preserving simple homotopy equivalence.
\end{lem}

\begin{rem} \label{rem_S}
From the result of Waldhausen \cite{Wal}, the Whitehead group $Wh(G_K)$ is trivial for a knot group $G_K$.
Therefore homotopy equivalence between finite CW-complexes whose fundamental groups are isomorphic to $G_K$ for some $K$ is simple homotopy equivalence. 
\end{rem}

Let $F$ be a field with (possibly trivial) involution $f \mapsto \bar{f}$.
We extend the involution to $F(t)$ by taking $t \mapsto t^{-1}$.
We equip $F(t)^{\oplus n}$ with the standard hermitian inner product $(\cdot, \cdot)$ defined by
\[ (v, w) := {^t v} \bar{w} \]
for $v,w \in F(t)^{\oplus n}$, where $^t v$ is the transpose of $v$.
For a representation $\rho \colon G_K \to GL_n(F)$, we define a representation $\rho^{\dagger} \colon G_K \to GL_n(F)$ by
\[ \rho^{\dagger}(\gamma) := \rho(\gamma^{-1})^{\ast} \]
for $\gamma \in G_K$, where $A^{\ast} := {^t \overline{A}}$ for a matrix $A$.

We can also refine the duality theorem for sign-determined Reidemeister torsion as follows.
\begin{thm} \label{thm_RD}
For a representation $\rho \colon G_K \to GL_n(F)$, if $H_*(E_K;F(t)_{\alpha \otimes \rho}^{\oplus n})$ vanishes, then so does $H_*(E_K;F(t)_{\alpha \otimes \rho^{\dagger}}^{\oplus n})$, and
\[ \widetilde{T}_{K, \rho^{\dagger}}(t) = (-1)^n \overline{\widetilde{T}_{K, \rho}(t)}. \]
\end{thm}

The proof is based on the following observations.
Let $(E'_K, \{ e'_i \})$ denote the PL manifold $E_K$ with the dual cell structure and choose a lift $\tilde{e}'_i$ which is the dual of $\tilde{e}_i$.
In the remainder of this section, for abbreviation, we write
\begin{align*}
C_q &:= C_q(\widetilde{E}_K) \otimes_{\alpha} \Q(t), &C_{\rho, q} &:= C_q(\widetilde{E}_K) \otimes_{\alpha \otimes \rho} F(t)^{\oplus n}, \\
C'_q &:= C_q(\widetilde{\partial E}_K) \otimes_{\alpha} \Q(t), &C'_{\rho, q} &:= C_q(\widetilde{\partial E}_K) \otimes_{\alpha \otimes \rho} F(t)^{\oplus n}, \\
C''_q &:= C_q(\widetilde{E}_K, \widetilde{\partial E}_K) \otimes_{\alpha} \Q(t), &C''_{\rho, q} &:= C_q(\widetilde{E}_K, \widetilde{\partial E}_K) \otimes_{\alpha \otimes \rho} F(t)^{\oplus n}, \\
D_q &:= C_q(\widetilde{E_K'}) \otimes_{\alpha} \Q(t), &D_{\rho, q} &:= C_q(\widetilde{E_K'}) \otimes_{\alpha \otimes \rho^{\dagger}} F(t)^{\oplus n}, \\
B'_q &:= \im(\partial : C'_{q+1} \to C'_q), &B'_{\rho, q} &:= \im(\partial : C'_{\rho, q+1} \to C'_{\rho, q}), \\
B''_q &:= \im(\partial : C''_{q+1} \to C'_q), &B''_{\rho, q} &:= \im(\partial : C''_{\rho, q+1} \to C''_{\rho, q}).
\end{align*}
Note that since a direct computation implies
\begin{equation} \label{eq_A}
H_*(\partial E_K; F(t)_{\alpha \otimes \rho}^{\oplus n}) = 0,
\end{equation}
we have
\begin{equation}
\begin{split}
\dim B'_{\rho, i} &= \sum_{j=0}^i (-1)^{i-j} \dim C'_{\rho, j} \\
&= \sum_{j=0}^i (-1)^{i-j} n \dim C'_j = n \dim B'_i.
\end{split} \label{eq_d1}
\end{equation}
(See, for example, \cite[Subsection 3.3.]{KL}.)
Similarly, if $H_*(E_K; F(t)_{\alpha \otimes \rho}^{\oplus n}) = 0$, then it follows from \eqref{eq_A} and the long exact sequence of the pair $(E_K, \partial E_K)$ that $H_*(E_K, \partial E_K; F(t)_{\alpha \otimes \rho}^{\oplus n}) = 0$, and so
\begin{equation}
\dim B''_{\rho, i} = n \dim B''_i. \label{eq_d2}
\end{equation}
The inner product
\[ [\cdot, \cdot] \colon C_q(\widetilde{E_K'}) \times C_{3-q}(\widetilde{E}_K, \widetilde{\partial E}_K) \to \Z[G_K] \]
defined by
\[ [\tilde{e}'_i, \tilde{e}_j] := \sum_{\gamma \in G_K} (\tilde{e}'_i, \tilde{e}_j \cdot \gamma^{-1}) \gamma, \]
where $(\cdot, \cdot)$ denote the intersection pairing, induces an inner product
\[ \langle \cdot, \cdot \rangle \colon D_{\rho, q} \times  C''_{\rho, 3-q} \to \C(t) \]
defined by
\[ \langle \tilde{e}'_i \otimes v, \tilde{e}_j \otimes w \rangle := (v, [\tilde{e}'_i, \tilde{e}_j] \cdot w) \]   
for $v, w \in \C(t)^{\oplus n}$.
(See, for example, \cite[Lemma 2.]{Mi}.)
This gives
\begin{equation}
D_{\rho, q} \cong (C''_{\rho, 3-q})^{\ast}. \label{eq_d3}
\end{equation}
The differential $\partial_q$ of $D_{\rho, q}$ corresponds with $(-1)^{q} \partial_{3-q}^{\ast}$ of $(C''_{\rho, 3-q})^{\ast}$ under this isomorphism. 
We also have
\begin{equation}
D_q \cong (C''_{3-q})^{\ast}. \label{eq_d4}
\end{equation}

\begin{lem} \label{lem_D}
For a representation $\rho \colon G_K \to GL_n(F)$,
\[ H_q(E_K; F(t)_{\alpha \otimes \rho^{\dagger}}^{\oplus n}) \cong H_{3-q}(E_K; F(t)_{\alpha \otimes \rho}^{\oplus n})^{\ast}. \]
\end{lem}

\begin{proof}
From \eqref{eq_d3} and the universal coefficient theorem,
\[ H_q(E_K; F(t)_{\alpha \otimes \rho^{\dagger}}^{\oplus n}) \cong H_{3 - q}(E_K, \partial E_K; F(t)_{\alpha \otimes \rho}^{\oplus n})^{\ast}. \]
From \eqref{eq_A} and the long exact sequence of the pair $(E_K, \partial E_K)$,
\[ H_*(E_K; F(t)_{\alpha \otimes \rho}^{\oplus n}) \cong H_*(E_K, \partial E_K; F(t)_{\alpha \otimes \rho}^{\oplus n}). \]
These isomorphisms prove the lemma.
\end{proof}

Now we prove the theorem.
\begin{proof}[Proof of Theorem \ref{thm_RD}]
Lemma \ref{lem_D} proves the first assertion.

In the following we use the notation in Definition \ref{def_NR}.
We choose an orthonormal basis $\langle f_1, \dots, f_n \rangle$ of $F(t)^{\oplus n}$ with respect to the hermitian product $(\cdot, \cdot)$ defined above.
Let $\boldsymbol{c}'$, $\boldsymbol{c}''$, $\boldsymbol{c}'_0$, $\boldsymbol{c}''_0$, $\tilde{\boldsymbol{c}}'$ and $\tilde{\boldsymbol{c}}''$ be the bases of $C_*(\partial E_K)$, $C_*(E_K, \partial E_K)$, $C'_*$, $C''_*$, $C'_{\rho, *}$ and $C''_{\rho, *}$ respectively induced by $\boldsymbol{c}$, $\tilde{\boldsymbol{c}}_0$ and $\tilde{\boldsymbol{c}}$.
We set
\begin{align*}
\boldsymbol{c}^{\ast} &:= \langle e'_1, \dots , e'_{dim C_*} \rangle, \\
\tilde{\boldsymbol{c}}_0^{\ast} &:= \langle e'_1 \otimes 1, \dots , e'_{dim C_*} \otimes 1 \rangle, \\
\tilde{\boldsymbol{c}}^{\ast} &:= \langle \tilde{e}'_1 \otimes f_1, \dots , \tilde{e}'_1 \otimes f_n, \dots , \tilde{e}'_{dim C_*} \otimes f_1, \dots , \tilde{e}'_{dim C_*} \otimes f_n \rangle.
\end{align*}
From \eqref{eq_d3} and the duality for algebraic torsion (\cite[Theorem 1.9]{T2}),
\[ \tor(D_{\rho, *}, \tilde{\boldsymbol{c}}^{\ast}) = (-1)^{\sum_i \dim B''_{\rho, i-1} \dim B''_{\rho, i}} \overline{\tor(C''_{\rho, *}, \tilde{\boldsymbol{c}}'')}. \]
On the other hand, from the exact sequence
\[ 0 \to C'_{\rho, *} \to C_{\rho, \ast} \to C''_{\rho, *} \to 0 \] 
and the multiplicativity for algebraic torsion (\cite[Theorem 1.5]{T2}),
\[ \tor(C_{\rho, *}, \tilde{\boldsymbol{c}}) = (-1)^{\sum_i \dim B'_{\rho, i-1} \dim B''_{\rho, i}} \tor(C'_{\rho, *}, \tilde{\boldsymbol{c}}') \tor(C''_{\rho, *}, \tilde{\boldsymbol{c}}''). \]
Therefore
\begin{equation}
\tor(C_{\rho, *}, \tilde{\boldsymbol{c}}) = (-1)^{\sum_i (\dim B'_{\rho, i-1} + \dim B''_{\rho, i-1})\dim B''_{\rho, i}} \tor(C'_{\rho, *}, \tilde{\boldsymbol{c}}') \overline{\tor(D_{\rho, *}, \tilde{\boldsymbol{c}}^{\ast})}. \label{eq_d5}
\end{equation}
Similarly,
\begin{equation}
\tor(C_*, \tilde{\boldsymbol{c}}_0) = (-1)^{\sum_i (\dim B'_{i-1} + \dim B''_{i-1})\dim B''_i} \tor(C'_*, \tilde{\boldsymbol{c}}'_0) \overline{\tor(D_*, \tilde{\boldsymbol{c}}_0^{\ast})}. \label{eq_d6}
\end{equation}
We set
\begin{align*}
d'' &:= \frac{1}{2}(\hdeg \tor(C'_*, \tilde{\boldsymbol{c}}'_0) + \ldeg \tor(C'_*, \tilde{\boldsymbol{c}}'_0)), \\
d^{\ast} &:= \frac{1}{2}(\hdeg \tor(D_*, \tilde{\boldsymbol{c}}_0^{\ast}) + \ldeg \tor(D_*, \tilde{\boldsymbol{c}}_0^{\ast})).
\end{align*}
From \eqref{eq_d6},
\begin{equation}
d' = d'' - d^{\ast}. \label{eq_d7}
\end{equation}
Since it is well-known that
\[ (t - 1) \tor (C_*, \tilde{\boldsymbol{c}}_0) \equiv \Delta_K(t) \mod \langle \pm t \rangle, \]
from Lemma \ref{lem_NR},
\[ \lim_{t \to 1} \tau_0 (t - 1) \tor (C_*, \tilde{\boldsymbol{c}}_0) = \lim_{t \to 1} \tau_0^{\ast} (t - 1) \tor (D_*, \tilde{\boldsymbol{c}}_0^{\ast}) = \pm1, \]
where
\[ \tau_0^{\ast} := \sgn \tor (C_*(E'_K; \R),\boldsymbol{c}^{\ast} ,\boldsymbol{h}). \]
Hence, by multiply \eqref{eq_d6} by $(t - 1)$ and taking $t \to 1$, we obtain
\begin{equation}
\tau_0 = -(-1)^{\sum_i (\dim B'_{i-1} + \dim B''_{i-1})\dim B''_i} \tau'_0 \tau_0^{\ast}, \label{eq_d8}
\end{equation}
where
\[ \tau'_0 := \lim_{t \to 1} \tor(C'_*, \tilde{\boldsymbol{c}}'_0). \]

From \eqref{eq_d1}, \eqref{eq_d2}, \eqref{eq_d5}, \eqref{eq_d7} and \eqref{eq_d8},
\begin{align*}
\widetilde{T}_{K, \rho}(t) &= \frac{\tau_0^n}{(\epsilon t^n)^{d'}} \tor(C_{\rho, *}, \tilde{\boldsymbol{c}}) \\
&= (-1)^n \frac{(\tau'_0)^n}{(\epsilon t^n)^{d''}} \tor(C'_{\rho, *}, \tilde{\boldsymbol{c}}') \cdot \overline{\frac{(\tau_0^{\ast})^n}{(\epsilon t^n)^{d^{\ast}}} \tor(D_{\rho, *}, \tilde{\boldsymbol{c}}^{\ast})}.
\end{align*}
A direct computation implies
\[ \tor(C'_*, \tilde{\boldsymbol{c}}'_0) = \tau'_0 t^{d''}. \]
(See, for example, \cite[Subsection 3.3.]{KL}.)
Since the normalized invariants do not change by conjugation of representations, we can assume that elements of $\rho(\pi_1 \partial E_K)$ are all diagonal.
This deduces
\[ \tor(C'_{\rho, *}, \tilde{\boldsymbol{c}}') = (\tau'_0)^n (\epsilon t^n)^{d''}. \] 
Thus
\[ \frac{(\tau'_0)^n}{(\epsilon t^n)^{d''}} \tor(C'_{\rho, *}, \tilde{\boldsymbol{c}}') = 1. \]
Further it can be easily seen that
\[ \frac{(\tau_0^{\ast})^n}{(\epsilon t^n)^{d^{\ast}}} \tor(D_{\rho, *}, \tilde{\boldsymbol{c}}^{\ast}) = \widetilde{T}_{K, \rho^{\dagger}}(t), \]
and the proof is complete.
\end{proof}

In the normalized setting, Theorem \ref{thm_AR} also holds.
\begin{thm} \label{thm_NAR}
For a representation $\rho \colon G_K \to GL_n(F)$ such that $H_*(E_K; F(t)_{\alpha \otimes \rho}^{\oplus n})$ vanishes,
\[ \widetilde{\Delta}_{K,\rho}(t) = \widetilde{T}_{K,\rho}(t). \]
\end{thm}

\begin{proof}
We choose a Wirtinger presentation $G_K = \langle x_1, \dots , x_m \mid r_1, \dots , r_{m-1} \rangle$ and take the CW-complex $W$ corresponding with the presentation.
Namely, $W$ has one vertex, $m$ edges labeled by the generators $x_1, \dots, x_m$ and $(m-1)$ 2-cells attached along the relations $r_1, \dots ,r_{m-1}$.
Let $x_1, \dots , x_m$ and $r_1, \dots ,r_{m-1}$ also denote the cells.
It is easy to see that $W$ is homotopy equivalent to $E_K$.
It follows from Remark \ref{rem_S} that $W$ is simple homotopy equivalent to $E_K$.
Thus from Lemma \ref{lem_NR} we can compute the normalized torsion $\widetilde{T}_{K, \rho}$ as that of $W$.

The chain complex $C_*(W;\R)$ is written as:
\[ 0 ~\to~ \bigoplus\limits_{j = 1}^{m-1} \R r_j ~\xrightarrow{\partial_2}~ \bigoplus\limits_{i = 1}^m \R x_i ~\xrightarrow{\partial_1}~ \R pt ~\to~ 0, \]
where
\begin{align*}
\partial_1 &= \boldsymbol{0}, \\
\partial_2 &= \left( \varphi \left( \frac{\partial r_j}{\partial x_i} \right) \right).
\end{align*}
Let $c_0 = pt$, $c_1 = \langle x_1, \dots , x_m \rangle$ and $c_2 = \langle r_1, \dots , r_{m-1} \rangle$.
We choose $b_1 = \partial c_2$ and $h_0 = [pt]$, $h_1 = [x_k]$ ($1 \leq k \leq m$).
Then
\begin{align*}
\tau_0 &= \sgn (-1)^{|C_*(W; \R)|} \frac{[b_1 \cup \tilde{h}_1 / c_1]}{[\tilde{h}_0 / c_0][\tilde{b}_1 / c_2]} \\
&= - \sgn \det
\begin{pmatrix}
& 0 \\
& \vdots \\
& 0 \\
\left( \varphi \left( \frac{\partial r_j}{\partial x_i} \right) \right) & 1 \\
& 0 \\
& \vdots \\
& 0
\end{pmatrix} \\
&= (-1)^{k + m + 1} \delta.
\end{align*}

We define an involution $\bar{\cdot} \colon \Z[G_K] \to \Z[G_K]$ by extending the inverse operation $\gamma \mapsto \gamma^{-1}$ of $G_K$ linearly.  
We can choose lifts $\widetilde{pt}$, $\tilde{x}_i$ and $\tilde{r}_j$ so that $C_*(\widetilde{W}) \otimes_{\alpha \otimes \rho} F(t)^{\oplus n}$ is written as:
\[ 0 ~\to~ \bigoplus\limits_{1 \leq j \leq m-1 , 1 \leq l \leq n} F(t) (\tilde{r}_j \otimes f_l) ~\xrightarrow{\tilde{\partial}_2}~ \bigoplus\limits_{1 \leq i \leq m , 1 \leq l \leq n}  F(t) (\tilde{x}_i \otimes f_l) ~\xrightarrow{\tilde{\partial}_1}~ \bigoplus\limits_{1 \leq l \leq n} F(t) (\widetilde{pt} \otimes f_l) ~\to~ 0, \]
where
\begin{align*}
\tilde{\partial}_1 (\tilde{x}_i \otimes f_l) &= \widetilde{pt} \otimes \Phi (\overline{\tilde{x}_i - 1}) f_l \\
\tilde{\partial}_2 (\tilde{r}_j \otimes f_l) &= \sum_{i=1}^m \tilde{x}_i \otimes \Phi \left( \overline{\frac{\partial r_j}{\partial x_i}} \right) f_l.
\end{align*}
Let $c_0' = \langle \widetilde{pt} \otimes f_1, \dots , \widetilde{pt} \otimes f_{n} \rangle $, $c_1' = \langle \tilde{x}_1 \otimes f_1, \dots ,\tilde{x}_1 \otimes f_n, \dots , \tilde{x}_m \otimes f_1, \dots , \tilde{x}_m \otimes f_n \rangle$ and $c_2' = \langle \tilde{r}_1 \otimes f_1, \dots ,\tilde{r}_1 \otimes f_n, \dots , \tilde{r}_{m-1} \otimes f_1, \dots , \tilde{r}_{m-1} \otimes f_n \rangle$.
We choose $b_0' = \partial \langle \tilde{x}_k \otimes f_1, \dots ,\tilde{x}_k \otimes f_n \rangle$ and $b_1' = \partial c_2'$.
Since $H_*(W;F(t)_{\alpha \otimes \rho}^{\oplus n})$ vanishes, $|C_*(\widetilde{W}) \otimes_{\alpha \otimes \rho} F(t)^{\oplus n}| = 0$, and so
\begin{align*}
\quad \tor (C_*(\widetilde{W}) \otimes_{\alpha \otimes \rho} F(t)^{\oplus n}, \langle \tilde{c}_0, \tilde{c}_1, \tilde{c}_2 \rangle ) &= \frac{[b_1' \cup \tilde{b}_0' / c_1']}{[b_0' / c_0'][\tilde{b}_1' / c_2']} \\
&= \frac{\det
\begin{pmatrix}
& 0 \\
& \vdots \\
& 0 \\
\left( \Phi \left( \overline{\frac{\partial r_j}{\partial x_i}} \right) \right) & I \\
& 0 \\
& \vdots \\
& 0
\end{pmatrix}
}{\det \Phi(\overline{x_k - 1})} \\
&= (-1)^{n(k+m)} \frac{\det \left( {^t \Phi \left( \overline{\frac{\partial r_i}{\partial x_j}} \right)} \right)}{\det {^t \Phi(\overline{x_k - 1})}}.
\end{align*}
Similarly, we obtain
\[ \tor (C_*(\widetilde{W}) \otimes_{\alpha} \Q (t), \langle \tilde{c}_0, \tilde{c}_1, \tilde{c}_2 \rangle) = (-1)^{(k+m)} \frac{\det \left( \tilde{\alpha} \left( \overline{\frac{\partial r_i}{\partial x_j}} \right) \right)}{\det \tilde{\alpha} (\overline{x_k - 1})}. \]
Hence $d' = -d$.

The above computations imply
\[ \widetilde{T}_{K, \rho}(t) = (-1)^n \overline{\widetilde{\Delta}_{K, \rho^{\dagger}}(t)}, \]
where we consider the trivial involution on $F$.
Now establishes the theorem follows from Theorem \ref{thm_RD}.
\end{proof}

From the above theorems and the following lemma, we have the duality theorem for normalized twisted Alexander invariants.
\begin{lem}
If $H_*(E_K; F(t)_{\alpha \oplus \rho}^{\oplus n})$ does not vanish, then
\[ \widetilde{\Delta}_{K, \rho}(t) = \widetilde{\Delta}_{K, \rho^{\dagger}}(t) = 0. \]
\end{lem}

\begin{proof}
If $H_*(E_K; F(t)_{\alpha \otimes \rho}^{\oplus n})$ does not vanish, then neither does $H_*(E_K; F(t)_{\alpha \otimes \rho^{\dagger}}^{\oplus n})$ from Lemma \ref{lem_D}.
Since
\[ \sum_{q=0}^2 \dim H_q(E_K; F(t)_{\alpha \oplus \rho}^{\oplus n}) = n \chi(E_K) = 0, \] 
it follows from the assumption and \eqref{eq_A} that $H_1(E_K; F(t)_{\alpha \otimes \rho}^{\oplus n}) \neq 0$, and so $\widetilde{\Delta}_{K, \rho}(t) = 0.$
Similarly, we obtain $\widetilde{\Delta}_{K, \rho^{\dagger}}(t) = 0$, which proves the lemma.
\end{proof}

\begin{thm} \label{thm_D}
For a representation $\rho \colon G_K \to GL_n(F)$,
\[ \widetilde{\Delta}_{K, \rho^{\dagger}}(t) = (-1)^n \overline{\widetilde{\Delta}_{K, \rho}(t)}. \]
\end{thm}

For a unitary representation $\rho$, the difference between the highest and lowest coefficients of $\Delta_{K, \rho}(t)$ is not clear from Theorem \ref{thm_UD} because of the ambiguity.
However, this difference is now strictly determined from the following corollary.
\begin{cor}
For a representation $\rho \colon G_K \to U(n) \text{ or } O(n)$,
\[ \widetilde{\Delta}_{K, \rho}(t) = (-1)^n \overline{\widetilde{\Delta}_{K, \rho}(t)}. \]
\end{cor}

\begin{exmp} \label{ex_T}
Let $K$ be the $(p,q)$ torus knot where $p,q > 1$ and $(p,q) = 1$.
It is well known that the knot group has a presentation
\[ G_K = \langle x, y \mid x^p y^{-q} \rangle \]
where $\hdeg \alpha (x) = q$ and $\hdeg \alpha (y) = p$.
The 2-dimensional complex
$W$ corresponding with this presentation is $K(G_K, 1)$.
Therefore we can use this presentation for the computation via Lemma \ref{lem_NR}, Remark \ref{rem_S} and Theorem \ref{thm_NAR}.

From the result of Klassen \cite{Kl}, all the irreducible $SU(2)$-representations up to conjugation are given as follows:
\begin{align*}
\rho_{a,b,s} \colon G_K &\to SU(2) : \\
x &\mapsto 
\begin{pmatrix}
\cos \frac{a \pi}{p} + i \sin \frac{a \pi}{p} & 0 \\
0 & \cos \frac{a \pi}{p} - i \sin \frac{a \pi}{p}
\end{pmatrix}, \\
y &\mapsto 
\begin{pmatrix}
\cos \frac{b \pi}{q} + i \sin \frac{b \pi}{q} \cos \pi s & \sin \frac{b \pi}{q} \sin \pi s \\
- \sin \frac{b \pi}{q} \sin \pi s & \cos \frac{b \pi}{q} - i \sin \frac{b \pi}{q} \cos \pi s
\end{pmatrix},
\end{align*}
where $a,b \in \N$, $1 \leq a \leq p-1$, $1 \leq b \leq q-1$, $a \equiv b \bmod 2$ and $0 < s < 1$.
The normalized twisted Alexander invariants associated to these representations are computed as follows:
\[ \widetilde{\Delta}_{K,\rho_{a,b,s}}(t) = \frac{(t^{\frac{pq}{2}} - (-1)^a t^{-\frac{pq}{2}})^2}{(t^p - 2 \cos \frac{b \pi}{q} + t^{-p})(t^q - 2 \cos \frac{a \pi}{p} + t^{-q})}. \]
\end{exmp}

%%%%%%% Section 6 %%%%%%%%%%%%%%%%%%%%%%%%%%%%%%%%%%%%%%%%%%%%%%%%%%%
\section{Applications}

Now we consider applications of the normalized invariants.
First we extend the result of Goda-Kitano-Morifuji and Friedl-Kim.
We denote by $g(K)$ the genus of $K$.

Their results are as follows.
\begin{thm}[\cite{GKM}]
For a fibered knot $K$ and a unimodular representation $\rho \colon G_K \to SL_{2n}(F)$, $\mathrm{c}(\Delta_{K, \rho})$ is well-defined and equals $1$.
\end{thm}
\begin{thm}[\cite{C},\cite{FK}] \label{thm_FK}
For a fibered knot $K$ and a representation $\rho \colon G_K \to GL_n(R)$, $\Delta^1_{K, \rho}$ is monic and $\deg \Delta_{K, \rho} = n(2 g(K) - 1)$, where a polynomial is said to be monic if both the highest and lowest coefficients are units.
\end{thm}

In the normalized setting, we have the following theorem.
\begin{thm} \label{thm_F}
For a fibered knot $K$ and a representation $\rho \colon G_K \to GL_n(R)$,
\begin{align*}
\deg \widetilde{\Delta}_{K, \rho} &= 2 \hdeg \widetilde{\Delta}_{K, \rho} = n(2 g(K) - 1), \\
\mathrm{c}(\widetilde{\Delta}_{K, \rho}) &= \mathrm{c}(\nabla_K)^n {\epsilon}^{g(K) - \frac{1}{2}}.
\end{align*}
\end{thm}

\begin{proof}
The equality $\deg \widetilde{\Delta}_{K, \rho} = n(2g(K)-1)$ can be obtained from Theorem \ref{thm_FK}.
Since we have $\widetilde{\Delta}_{K,\iota \circ \rho} = \widetilde{\Delta}_{K, \rho}$, where $\iota$ is the natural inclusion $GL_n(R) \hookrightarrow GL_n(Q(R))$, we can assume $R$ is a field $F$.

Let $\psi$ denote the automorphism of a surface group induced by the monodromy map.
We can take the following presentation of $G_K$ by using the fibered structure:
\[ G_K = \langle x_1, \dots , x_{2g}, h \mid r_i := hx_ih^{-1} \psi_*(x_i)^{-1}, 1 \leq i \leq 2 g(K) \rangle \]
where $\alpha(x_i) = 1$ for all $i$ and $\alpha(h) = t$.
It is easy to see that the corresponding CW-complex is homotopy equivalent to the exterior $E_K$.
Thus we can compute the invariant by using the presentation as in Example \ref{ex_T}.

Since
\begin{align*}
\frac{\partial r_i}{\partial x_j} = 
\begin{cases}
h - \frac{\partial \psi_*(x_i)}{\partial x_i} &i = j \\
- \frac{\partial \psi_*(x_i)}{\partial x_j} &i \neq j
\end{cases},
\end{align*}
we have
\begin{align*}
\det A_{\tilde{\alpha}, 2g(K)+1} &= t^{2g(K)} + \dots + 1, \\
\det A_{\Phi,2g+1} &= \epsilon^{2g(K)} t^{2ng(K)} + \dots + (-1)^n \det (\Phi (\frac{\partial \psi_*(x_i)}{\partial x_j})), \\
\det \Phi (h - 1) &= \epsilon t^n + \dots + (-1)^n .
\end{align*}
From the classical theorem of Neuwirth, which states that the degree of the Alexander polynomial of a fibered knot equals twice the genus, we can determine that the lowest degree term of the first equality equals $1$. 
Since 
\begin{align*}
\delta &= \sgn \left. \mathrm{c}(\nabla_K)\nabla_K(t^{\frac{1}{2}}-t^{-\frac{1}{2}}) \right|_{t=1} \\
&= \mathrm{c}(\nabla_K) \\
d &= g(K) - \frac{1}{2},
\end{align*}
we obtain $\hdeg \widetilde{\Delta}_{K, \rho} = n(g(K) - \frac{1}{2})$ and $\mathrm{c}(\widetilde{\Delta}_{K, \rho}) = \mathrm{c}(\nabla_K)^n {\epsilon}^{2 g(K) - 1}$.
\end{proof}

Next we study behavior of the highest degrees of the normalized invariants.
\begin{defn}
A Seifert surface for a knot $K$ is said to be \textit{canonical} if it is obtained from a diagram of $K$ by applying the Seifert algorithm.
The minimum genus over all canonical Seifert surfaces is called the \textit{canonical genus} and denoted by $g_c(K)$.
A Seifert surface $S$ is said to be \textit{free} if $\pi_1(S^3 \setminus S)$ is a free group.
This condition is equivalent to that $S^3 \setminus N(S)$ is a handlebody, where $N(S)$ is an open regular neighborhood of $S$.
The minimum genus over all free Seifert surfaces is called the \textit{free genus} and denoted by $g_f(K)$.
\end{defn}

\begin{rem}
Since every canonical Seifert surface is free, the following fundamental inequality holds:
\[ g(K) \leq g_f(K) \leq g_c(K). \]
\end{rem}

The highest degrees of the normalized invariants give lower bounds on the free genus.
\begin{thm} \label{thm_H}
For a representation $\rho \colon G_K \to GL_n(R)$, the following inequality holds:
\[ 2 \hdeg \widetilde{\Delta}_{K, \rho} \leq n (2 g_f(K) - 1). \]
\end{thm}

\begin{proof}
We choose a free Seifert surface $S$ with genus $g_f(K)$ and take a bicollar $S \times [-1, 1]$ of $S$ such that $S \times 0 = S$.
Let $\iota_{\pm} \colon S \hookrightarrow S^3 \setminus S$ be the embeddings whose images are $S \times \{ \pm 1 \}$.
Picking generator sets $\{ a_1, \dots, a_{2 g_f(K)} \}$ of $\pi_1 S$ and $\{ x_1, \dots, x_{2 g_f(K)} \}$ of $\pi_1(S^3 \setminus S)$ and setting $u_i := (\iota_+)_*(a_i)$ and $v_i := (\iota_-)_*(a_i)$ for all $i$, we have a presentation
\[ \langle x_1, \dots , x_{2g_f(K)}, h \mid r_i := hu_ih^{-1}v_i^{-1}, 1 \leq i \leq 2g_f(K) \rangle \]
of $G_K$ where $\alpha(x_i) = 1$ for all $i$ and $\alpha(h) = t$.

Collapsing surfaces $S \times *$ and the handlebody $S^3 \setminus (S \times [-1, 1])$ to bouquets, we can realize the $2$-dimensional complex corresponding with this presentation as a deformation retract of $E_K$.
Therefore we can compute the invariant by using the presentation as in Example \ref{ex_T}.
Since
\[ \frac{\partial r_i}{\partial x_j} = h \frac{\partial u_i}{\partial x_j} - \frac{\partial v_i}{\partial x_j}, \]
we have
\begin{align*}
\hdeg \widetilde{\Delta}_{K,\rho} &= \hdeg \det A_{\Phi,2g_f(K) + 1} - nd - n \\
&\leq 2 n g_f(K) - nd - n.
\end{align*}
Thus the proof is completed by showing that $d = g_f(K) - \frac{1}{2}$.

Let $V$ be the Seifert matrix with respect to the basis $\langle [a_1], \dots, [a_{2g_f(K)}] \rangle$ of $H_1(S; \Z)$ and $\langle [a_1]^*, \dots, [a_{2g_f(K)}]^* \rangle$ the dual basis of  $H_1(S^3 \setminus S; \Z)$, i.e., $lk([a_i], [a_j]^*) = \delta_{i, j}$.
We denote by $A_{\pm}$ the matrices representing $(\iota_{\pm})_* \colon H_1(S; \Z) \to H_1(S^3 \setminus S; \Z)$ with respect to the bases $\langle [a_1], \dots, [a_{2g_f(K)}] \rangle$ and $\langle [x_1], \dots, [x_{2g_f(K)}] \rangle$ and by $P$ the base change matrix of $H_1(S^3 \setminus S; \Z)$ from $\langle [x_1], \dots, [x_{2g_f(K)}] \rangle$ to $\langle [a_1]^*, \dots, [a_{2g_f(K)}]^* \rangle$.
It is well known that the matrices representing $(\iota_+)_*$ and $(\iota_-)_* \colon H_1(S; \Z) \to H_1(S^3 \setminus S; \Z)$ with respect to the bases $\langle [a_1], \dots, [a_{2g_f(K)}] \rangle$ and $\langle [a_1]^*, \dots, [a_{2g_f(K)}]^* \rangle$ are $V$ and $^t V$.
Hence
\begin{align*}
\det A_{\tilde{\alpha}, 2g_f(K) + 1} &= \det(t {^t A_+} - {^t A_-}) \\
&= \det(t A_+ - A_-) \\
&= \det(t PV - P {^t V}) \\
&= \pm \det(t V - {^t V}),
\end{align*}
and so $d = g_f(K) - \frac{1}{2}$ as required.
\end{proof}

\begin{exmp}
Let $K$ be the knot $11_{n73}$ illustrated in Figure~\ref{fig}.
The normalized Alexander polynomial of $K$ equals $t^2 - 2t + 3 - 2t^{-1} + t^{-2}$.

\begin{figure}[h]
\centering
\includegraphics[width=6cm]{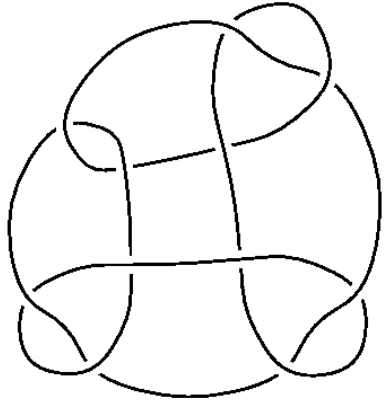}
\caption{The knot $11_{n73}$}
\label{fig}
\end{figure}

The Wirtinger presentation of the diagram in Figure~\ref{fig} consists of $11$ generators and $10$ relations:
\begin{align*}
&x_5 x_1 x_{5}^{-1} x_2^{-1}, &x_{11} x_2 x_{11}^{-1} x_3^{-1}, \\
&x_9 x_4 x_{9}^{-1} x_{3}^{-1}, &x_7 x_5 x_7^{-1} x_4^{-1}, \\
&x_1 x_5 x_1^{-1} x_6^{-1}, &x_8 x_7 x_8^{-1} x_6^{-1}, \\
&x_5 x_8 x_5^{-1} x_7^{-1}, &x_{10} x_9 x_{10}^{-1} x_8^{-1}, \\
&x_4 x_{10} x_4^{-1} x_9^{-1}, &x_2 x_{10} x_2^{-1} x_{11}^{-1}.
\end{align*}
Let $\rho \colon G_K \to SL_2(\F_2)$ be a nonabelian representation over $\F_2$ defined as follows:
\[ \rho(x_i) =
\begin{cases}
\begin{pmatrix}
1 & 0 \\
1 & 1 
\end{pmatrix}
~,~ \text{if } i=4,8 \\
\begin{pmatrix}
0 & 1 \\
1 & 0 
\end{pmatrix}
~,~ \text{if } i=7,9 \\
\begin{pmatrix}
1 & 1 \\
0 & 1 
\end{pmatrix}
~,~ \text{otherwise}\\ 
\end{cases}.
\]
From them, we can compute the normalized twisted Alexander invariant $\widetilde{\Delta}_{K,\rho}$ as:
\[ \widetilde{\Delta}_{K,\rho}(t) = t^5 + t + t^{-1} + t^{-5}. \]

Since $\deg \widetilde{\Delta}_{K,\rho} \neq 2 (\deg \Delta_K - 1)$, we can see that $K$ is not fibered.
Moreover, from Theorem \ref{thm_H}, we have
\[ 10 \leq 2 (2 g_f(K) - 1), \]
which becomes
\[ g_f(K) \geq 3. \]
On the other hand, we obtain a canonical Seifert surface with genus $3$ by applying the Seifert algorithm to the diagram in Figure~\ref{fig}.
Hence
\[ g_f(K) \leq g_c(K) \leq 3. \]
By these inequalities we conclude that
\[ g_f(K) = g_c(K) = 3. \]
\end{exmp}

\begin{rem}
Friedl and Kim \cite{FK} showed the following inequality:
\[ \deg \Delta_{K, \rho} \leq n(2 g(K) -1). \]
Therefore $g(K)$ also equals $3$ in the above example.
\end{rem}

%%%%%%% Reference %%%%%%%%%%%%%%%%%%%%%%%%%%%%%%%%%%%%%%%%%%%%%%%%%%%

\end{document}